\documentclass[10pt]{article}
\usepackage{amsmath,amssymb,amsthm}
\newcommand{\g}{\mathfrak{g}}
\newcommand{\sfrak}{\mathfrak{s}}
\newcommand{\Z}{\mathbb{Z}}

\newcommand{\bbhs}[6]{ {}_{#1}\psi_{#2} \left[ \genfrac{}{}{0pt}{}{#3}{#4} ; {#5},{#6} \right]}
\newcommand{\qbin}[3]{ \genfrac{[}{]}{0pt}{}{#1}{#2}_{#3} }

\numberwithin{equation}{section}
\allowdisplaybreaks

\newtheorem*{jtp}{Triple product identity}
\newtheorem*{qpi}{Quintuple product identity}
\newtheorem{theorem}{Theorem}

\title{A classical $q$-hypergeometric approach to the\\ $A_2^{(2)}$ standard modules}
\date{\today}
\author{Andrew V. Sills}

\begin{document}
\maketitle

\centerline{\em Dedicated to Krishna Alladi on the occasion of his sixtieth birthday}

\abstract{
This is a written expansion of the talk delivered by the author at the International Conference
on Number Theory in Honor of Krishna Alladi for his 60th Birthday, held at the University of
Florida, March 17--21, 2016.

  Here we derive Bailey pairs that give rise to Rogers--Ramanujan type identities which
are the principally specialized character of the $A_2^{(2)}$ standard
module $(\ell - 2i+2)\Lambda_0 + (i-1)\Lambda_1$ for any level $\ell$, and $i=1, 2$.}

\section{Notation and Motivation}
\subsection{$q$-series notation and classical results}  Let $q$ denote a formal variable.
The standard notation for the infinite rising $q$-factorial is
\[ (a;q)_\infty :=\prod_{j=0}^\infty (1-aq^j). \] 
In order to allow for positive and negative values of $n$, we define the finite rising 
$q$-factorial as
\[ (a;q)_n := \frac{(a;q)_\infty}{(aq^n;q)_\infty}. \]

We will also use the abbreviations $(q)_n$ and $(q)_\infty$ for $(q;q)_n$ and $(q;q)_\infty$
respectively.  Additionally,
\[ (a_1, a_2, \dots, a_r; q)_n := (a_1;q)_n (a_2;q)_n \cdots (a_r;q)_n \] and
\[ (a_1, a_2, \dots, a_r; q)_\infty := (a_1;q)_\infty (a_2;q)_\infty \cdots (a_r;q)_\infty. \]

The bilateral basic hypergeometric series is given by
\[ \bbhs{t}{t}{a_1,a_2,\dots,a_t}{b_1,b_2,\dots,b_t}{q}{z}  := \sum_{r\in\Z} 
\frac{(a_1;q)_r (a_2;q)_r \cdots (a_t;q)_r}{(b_1;q)_r (b_2;q)_r \cdots (b_t;q)_r}  z^r . \]

The $q$-binomial coefficient is
\[ \qbin{n}{m}{q} := \left\{ \begin{array}{ll} \frac{(q)_n}{(q)_m (q)_{n-m}}&
   \mbox{if $0\leq m\leq n $} \\ 0 & \mbox{otherwise} \end{array}
  \right. . \]

We will require the following classical results.  For our purposes, $z = \pm q^{r}$ for some
$r\in\frac 12 \mathbb{Z}$.
\begin{jtp}[Jacobi]  \cite[p. 15, Eq. (1.6.1)]{GR04}. 
 \begin{equation} \label{JTP}
  \sum_{n\in\Z} (-1)^n z^n q^{n^2} = (q/z, zq, q^2;q^2)_\infty. \end{equation}
\end{jtp}
\begin{qpi}[Fricke~\cite{F16}] cf.~\cite[p. 147, ex. 5.6]{GR04} 
\begin{multline} \label{QPI}
(-qz^3, -q^2 z^{-3}, q^3; q^3)_\infty - z (-qx^{-3}, -q^2 z^3, q^3; q^3)_\infty
\\= (q/z, z, q; q)_\infty (q/z^2, qz^2; q^2)_\infty. \end{multline}
\end{qpi}

\subsection{Certain affine Kac--Moody Lie algebras and their connection to $q$-series}
Let $\g$ denote the affine Kac--Moody Lie algebra $A_1^{(1)}$ or $A_2^{(2)}$ and
let $h_0, h_1$ denote the usual basis of a maximal toral subalgebra $T$ of $\g$.
Let $d$ denote the {degree derivation} of $\g$
and let $\tilde{T}:= T \oplus \mathbb C d$.
For all dominant integral $\lambda\in\tilde{T}^*$, there is a
unique irreducible, integrable, highest weight module $L(\lambda)$,
assuming (without loss of generality) that $\lambda(d) = 0$.
Also, $\lambda= s_0 \Lambda_0 + s_1 \Lambda_1$ where
$\Lambda_0$ and $\Lambda_1$ are the {fundamental weights},
given by $\Lambda_i(h_j) = \delta_{ij}$ and $\Lambda_i(d) = 0$;
$s_0$ and $s_1$ are nonnegative integers.
For $A_1^{(1)}$, the canonical central element is
$c= h_0 + h_1$, and for
$A_2^{(2)}$, the canonical central element is
$c = h_0 + 2h_1$.
The \emph{level} $\lambda(c)=\lambda_{\g}(c)$ of $L(\lambda)$ is
\[  \lambda(c) = \left\{   \begin{array}{ll} 
  s_0 + s_1 & \mbox{if $\g = A_1^{(1)},$} \\
  s_0 + 2s_1 &\mbox{if $\g = A_2^{(2)},$} \end{array}
 \right. \]
%
(cf. \cite{K90}, \cite{LM78}).
For brevity, it is common to refer to $L(\lambda) = L(s_0\Lambda_0 + s_1\Lambda_1)$ as
the ``$(s_0, s_1)$-module."

Additionally (\cite{LM78}),  there is an
infinite product $F_{\g}$ associated with $\g$, 
sometimes
called the {``fudge factor}," which needs to be divided out of the
the principally specialized character
$\chi(L(\lambda))
= \chi(s_0 \Lambda_0 + s_1\Lambda_1)$, in order to obtain
the quantities of interest here.  For $\g=A_1^{(1)}$,
the fudge
factor is given by
\[ F_{\g} = \left\{  \begin{array}{ll}(q;q^2)_\infty^{-1} &\mbox{if $\g = A_1^{(1)}$,} \\
\left[ (q;q^6)_\infty (q^5;q^6)_\infty \right]^{-1} &\mbox{if $\g = A_2^{(2)}$.  }
\end{array}
 \right. \]


  Also, $\g$ has a certain infinite-dimensional Heisenberg subalgebra
known as the 
{principal Heisenberg vacuum subalgebra} $\sfrak$
(consult~\cite{LW78} for the construction of $A_1^{(1)}$ and~\cite{KKLW81}
for that of $A_2^{(2)}$).
As demonstrated in~\cite{LW82}, the principal character
$\chi(\Omega(s_0 \Lambda_0 + s_1 \Lambda_1))$, where $\Omega(\lambda)$
is the vacuum space for $\sfrak$ in $L(\lambda)$, is
\begin{equation} \label{char}
\chi(\Omega(s_0 \Lambda_0 + s_1\Lambda_1))
=  \frac{\chi( L(s_0 \Lambda_0 + s_1 \Lambda_1)) }{ F_{\g}  },
\end{equation}
where $\chi(L(\lambda))$ is the principally specialized character of
$L(\lambda)$.

By~\cite{LM78} applied to~\eqref{char} in the case of $A_1^{(1)}$,
the standard modules of odd level 
correspond to 
Andrews' analytic generalization of the Rogers--Ramanujan
identities~\cite{A74}, known as the ``Andrews--Gordon identity," and the partition
theoretic generalization of the Rogers--Ramanujan identities due to B. Gordon~\cite{G61}.  Bressoud's even modulus counterpart to the Andrews--Gordon identity~\cite[p. 15, Eq. (3.4)]{B80}
and its partition theoretic counterpart~\cite[p. 64, Theorem, $j=0$ case]{B79};
was explained
vertex-operator theoretically in~\cite{LW84} and~\cite{LW85} to correspond to
the standard modules of even level in $A_1^{(1)}$.

The combined Andrews--Gordon--Bressoud identity (for both even and odd moduli) 
and its correspondence to the level $\ell$ standard modules of $A_1^{(1)}$ 
can be stated compactly as 
\begin{multline}\label{AndGorBress}
\chi(\Omega( (\ell+1-i)\Lambda_0 + (i-1)\Lambda_1 )) \\ = 
\sum_{n_1 \geqq n_2 \geqq \dots \geqq n_{k-1}\geqq 0}
 \frac{ q^{n_1^2 + n_2^2 + \cdots + n_{k-1}^2 + n_i+n_{i+1}+\cdots+n_{k-1}}}
 {(q)_{n_1-n_2} (q)_{n_2-n_3} \cdots (q)_{n_{k-2}-n_{k-1}}  (q)_{n_{k-1}} (-q)_{[[ 2\mid \ell ]] n_{k-1}} }
\\ = \frac{(q^i,q^{\ell+2-i},q^{\ell+2};q^{\ell+2})_\infty }{(q)_\infty },
\end{multline}
where $k:= k(\ell) =1 + \lfloor \ell/2 \rfloor$, $1\leqq i \leqq k$, and 
  \[ [[ P ]] := \left\{ \begin{array}{ll} 1 & \mbox{if $P$ is true,} \\
     0 & \mbox{if $P$ is false} \end{array} \right. . \]


A pair of sequences $(\alpha_n(a,q), \beta_n(a,q))$ form a \emph{Bailey pair} with respect to $a$ if
\[ \beta_n(a,q) = \sum_{s=0}^n \frac{\alpha_n(a,q)}{(q)_{n-s} (aq;q)_{n+s}} .\]
~\cite[p. 25--26]{A86}; cf.~\cite[pp. 2, 5]{B49}.

It is well known that identities of Rogers--Ramanujan type may be
derived by the insertion of Bailey pairs into limiting cases of 
Bailey's lemma~\cite[p. 25, Thm. 3.3; p. 27, Eq. (3.33)]{A86} such as
\begin{equation} \label{WBL}
\sum_{n=0}^\infty a^n q^{n^2} \beta_n(a,q) = \frac{1}{(aq;q)_\infty} \sum_{n=0}^\infty 
a^n q^{n^2} \alpha_n(a,q),
\end{equation}
and setting $a$ equal to a power of $q$.

   An efficient method for deriving~\eqref{AndGorBress} for odd $\ell$ is via
the Bailey lattice~\cite{AAB87}, which is an extension of the 
Bailey chain~(\cite{A84}; cf. \cite[\S 3.5, pp. 27ff]{A86})
built upon the ``unit Bailey pair"
\[ \beta_n(1,q) = \left\{
   \begin{array}{ll}
      1 &\mbox{if $n=0$}\\
      0 &\mbox{if $n>0$}
     \end{array}  \right. \]
\[ \alpha_n(1,q) = \left\{
     \begin{array}{ll}
        1 &\mbox{if $n=0$}\\
        (-1)^n q^{n(n-1)/2} (1+q^n) &\mbox{if $n>0$.}
     \end{array} \right. \]

  Similarly, for even $\ell$,
\eqref{AndGorBress} follows from a Bailey lattice built upon the
Bailey pair

\[ \beta_n(1,q) = \frac{1}{(q^2;q^2)_n}, \]
\[ \alpha_n(1,q) = \left\{
     \begin{array}{ll}
        1 &\mbox{if $n=0$}\\
        (-1)^n 2 q^{n^2} &\mbox{if $n>0$.}
     \end{array} \right. \]

Thus the standard modules of $A_1^{(1)}$ 
correspond to two interlaced instances of the Bailey lattice.

  In contrast, the standard modules of $A_{2}^{(2)}$ are not as well
understood, and a uniform $q$-series and partition correspondence
analogous to what is known for $A_1^{(1)}$ has to date remained
beyond our reach.

As with $A_1^{(1)}$, there are $1+\lfloor \frac{\ell}{2}
\rfloor$ inequivalent level $\ell$ standard modules associated with
the Lie algebra $A_2^{(2)}$, but the 
principal characters for the level $\ell$ standard modules
are given by instances of the quintuple product identity~\eqref{QPI}
(rather than the triple product identity) divided by $(q)_\infty$:

\begin{multline} \label{A22prodside}
\chi (\Omega( (\ell-2i+2)\Lambda_0 + (i-1)\Lambda_1 ))
\\ =  \frac{ (q^i, q^{\ell+3-i}, q^{\ell+3}; q^{\ell+3})_\infty (q^{\ell+3-2i},
  q^{\ell+2i+3}; q^{2\ell+6})_\infty}{(q)_\infty},
\end{multline}
where $1\leqq i \leqq 1 + \lfloor \frac{\ell}{2} \rfloor $; see~\cite{LM78}.

\section{Bailey pairs for $A_2^{(2)}$}

Let $( \alpha_n^{(\ell,i)}, \beta_n^{(\ell,i)})$ denote the Bailey pair which, upon insertion into
\eqref{WBL} with $a=1$, gives the principally specialized character of the $A_2^{(2)}$ standard
module  $(\ell-2i+2)\Lambda_0 + (i-1)\Lambda_1$.

\subsection{Bailey pairs for $\chi (\Omega( \ell\Lambda_0)) $}
 \[ \alpha_n = \alpha_n^{(\ell,1)}(1,q) = \left\{ 
\begin{array}{ll} 1 & \mbox{if $n=0$}\\
q^{\frac 32(\ell-3)r^2 - \frac12 (\ell-3) r} +q^{\frac 32(\ell-3)r^2 + \frac 12(\ell-3)r } & \mbox{if $n=3r>0$ } \\
-q^{\frac 32(\ell-3)r^2 + \frac12 (\ell-3) r} &\mbox{if $n=3r+1$} \\
-q^{\frac 32(\ell-3)r^2 - \frac12 (\ell-3) r} &\mbox{if $n=3r-1$}
\end{array}
 \right.
\]

\begin{align}
\beta_{n}^{(\ell,1)} (1,q) & = \sum_{s=0}^n \frac{\alpha_s^{(\ell,1)} (1,q) }{(q)_{n-s} (q)_{n+s}} \notag \\
                           & = \frac{\alpha_0}{(q)_n^2} 
                           + \sum_{r\geq 1}\frac{\alpha_{3r}}{(q)_{n-3r} (q)_{n+3r} }
                           + \sum_{r\geq 0}\frac{\alpha_{3r+1}}{(q)_{n-3r-1} (q)_{n+3r+1} } \notag \\ &\qquad 
                            + \sum_{r\geq 1}\frac{\alpha_{3r-1}}{(q)_{n-3r+1} (q)_{n+3r-1} } \notag\\
                     & = \frac{1}{(q)_n^2} 
                       + \sum_{r\geq 1}\frac{  q^{\frac 32(\ell-3)r^2 - \frac12 (\ell-3) r} +q^{\frac 32(\ell-3)r^2 + \frac 12(\ell-3)r } }{(q)_{n-3r} (q)_{n+3r} } \notag\\
                        & \qquad   + \sum_{r\geq 0}\frac{q^{\frac 32(\ell-3)r^2 + \frac12 (\ell-3) r}}{(q)_{n-3r-1} (q)_{n+3r+1} }
                            + \sum_{r\geq 1}\frac{ q^{\frac 32(\ell-3)r^2 - \frac12 (\ell-3) r} }{(q)_{n-3r+1} (q)_{n+3r-1} } \notag\\
                     &=\sum_{r\in\Z}\frac{  q^{\frac 32(\ell-3)r^2 + \frac 12(\ell-3)r } }{(q)_{n-3r} (q)_{n+3r} } 
                          -\sum_{r\in\Z}\frac{q^{\frac 32(\ell-3)r^2 + \frac12 (\ell-3) r}}{(q)_{n-3r-1} (q)_{n+3r+1} }
                          \notag \\
                  & =  \sum_{r\in\Z}\frac{  q^{\frac 32(\ell-3)r^2 + \frac 12(\ell-3)r } }{(q)_{n-3r} (q)_{n+3r+1} }
                  \Big(  (1-q^{n+3r+1}) - (1-q^{n-3r})  \Big)  \notag    \\
                  & =  \sum_{r\in\Z}\frac{  q^{\frac 32(\ell-3)r^2 + \frac 12(\ell-3)r } }{(q)_{n-3r} (q)_{n+3r+1} }
                ( q^{n-3r} - q^{n+3r+1}) \notag \\ 
                 & =  \sum_{r\in\Z}\frac{  q^{\frac 32(\ell-3)r^2 + \frac 12(\ell-3)r + n - 3r} }{(q)_{n-3r} (q)_{n+3r+1} }
                ( 1- q^{6r+1} ) \notag \\       
                 & = \frac{q^n}{ (q)_n (q)_{n+1} } \sum_{r\in\Z}\frac{  q^{\frac 32(\ell-3)r^2 + \frac 12(\ell-9)r }
                  ( 1- q^{6r+1} ) (q^{n-3r+1};q)_{3r}   }{ (q^{n+2};q)_{3r} }\notag \\
                  & =  \frac{q^n}{ (q)_n (q^2;q)_{n} } 
                  \sum_{r\in\Z}\frac{  
                  ( 1- q^{6r+1} ) (q^{-n};q)_{3r}   }{ (1-q) (q^{n+2};q)_{3r} } (-1)^r q^{\frac 32(\ell-6)r^2 + \frac 12(\ell-6)r } \label{VWPpsi}.
\end{align}

For each $\ell = 1,2,3,\dots$, the series expression in~\eqref{VWPpsi} is a limiting case of
a very-well-poised bilateral basic hypergeometric series.  

For example, we have 
\begin{multline*} \frac{q^{-n} (q)_n (q)_{n+1}}{1-q} \beta^{(\ell,1)}(1,q) \\= 
\left\{ \begin{array}{ll}
 \displaystyle{  \underset{e\to 0}{\lim} \bbhs{8}{8}
  {q^{\frac 72},-q^{\frac 72},q^{-n},q^{1-n},q^{2-n},e,e,e}
  {q^{\frac 12},-q^{\frac 12},q^{n+4},q^{n+3},q^{n+2},\frac{q^4}{e},\frac{q^4}{e},\frac{q^4}{e}}{q^3}
  {\frac{q^{3n+6}}{e^3} } }&\mbox{if $\ell=3$}\\
      \displaystyle{   \underset{e\to 0}{\lim} \bbhs{8}{8}
  {q^{\frac 72},-q^{\frac 72},q^{-n},q^{1-n},q^{2-n},e,e,-q^2}
  {q^{\frac 12},-q^{\frac 12},q^{n+4},q^{n+3},q^{n+2},\frac{q^4}{e},\frac{q^4}{e},-q^2}{q^3}
  { \frac{-q^{3n+4}}{e^2} } }&\mbox{if $\ell=4$}\\
      \displaystyle{  \underset{e\to 0}{\lim} \bbhs{6}{6}
  {q^{\frac 72},-q^{\frac 72},q^{-n},q^{1-n},q^{2-n},e}
  {q^{\frac 12},-q^{\frac 12},q^{n+4},q^{n+3},q^{n+2},\frac{q^4}{e}}{q^3}{ \frac{q^{3n+2}}{e} }
  } &\mbox{if $\ell=5$}\\
      \displaystyle{  \bbhs{6}{6}{q^{\frac 72},-q^{\frac 72},q^{-n},q^{1-n},q^{2-n},-q^2}
  {q^{\frac 12},-q^{\frac 12},q^{n+4},q^{n+3},q^{n+2},-q^2}{q^3}{ -q^{3n} } }&\mbox{if $\ell=6$}\\  
        \displaystyle{ \underset{e\to\infty}{\lim} \bbhs{6}{6}
  {q^{\frac 72},-q^{\frac 72},q^{-n},q^{1-n},q^{2-n},e}
  {q^{\frac 12},-q^{\frac 12},q^{n+4},q^{n+3},q^{n+2},\frac{q^4}{e}}{q^3}{  \frac{q^{3n+2}}{e}  } }&\mbox{if $\ell=7$} \\
       \displaystyle{ \underset{e\to \infty}{\lim} \bbhs{8}{8}
  {q^{\frac 72},-q^{\frac 72},q^{-n},q^{1-n},q^{2-n},e,e,-q^2}
  {q^{\frac 12},-q^{\frac 12},q^{n+4},q^{n+3},q^{n+2},\frac{q^4}{e},\frac{q^4}{e},-q^2}{q^3}
  { \frac{-q^{3n+4}}{e^2} } }&\mbox{if $\ell=8$}\\
        \displaystyle{  \underset{e\to \infty}{\lim} \bbhs{8}{8}
  {q^{\frac 72},-q^{\frac 72},q^{-n},q^{1-n},q^{2-n},e,e,e}
  {q^{\frac 12},-q^{\frac 12},q^{n+4},q^{n+3},q^{n+2},\frac{q^4}{e},\frac{q^4}{e},\frac{q^4}{e}}{q^3}
  {\frac{q^{3n+6}}{e^3} } }&\mbox{if $\ell=9$}\\
\end{array}.
\right. 
\end{multline*}

Observe that the easiest cases are $\ell = 5,6,7$, as these are instances of Bailey's summable
bilateral very-well-poised ${}_6\psi_6$~\cite[Eq. (4.7)]{B36}; cf.~\cite[p. 357, Eq. (II.33)]{GR04}.  Indeed, Slater evaluated the cases $\ell = 5$ and $7$~\cite[p. 464, 
Eqs. (3.4) and (3.3) resp.]{S51}, while McLaughlin and Sills evaluated the case $\ell=6$~\cite[p. 772, Table 3.1, line (P2)]{MS08}.

We have
\begin{align}
\beta^{(5,1)}_n (1,q) &= \frac{q^{n^2}}{(q)_{2n}}, \label{beta51}\\
\beta^{(6,1)}_n (1,q) &= \frac{ q^n (-1;q^3)_n }{ (-1;q)_n (q)_{2n}  }, \label{beta61}\\
\beta^{(7,1)}_n (1,q) &= \frac{q^n }{(q)_{2n}} . \label{beta71}
\end{align}

To evaluate $\beta^{(\ell,1)}_n (1,q)$ for levels $\ell=3,4,8,9$, we can use the following
identity~\cite[p. 147, exercise 5.11]{GR04}, analogous to Bailey's ${}_6\psi_6$ sum:  for nonnegative
integer $n$, 
\begin{multline} \bbhs{8}{8}{q\sqrt{a},-q\sqrt{a},c,d,e,f,a q^{-n}, q^{-n}}{ \sqrt{a},-\sqrt{a}, aq/c, aq/d, aq/e,aq/f, q^{n+1},
a q^{n+1} }{q}{\frac{a^2 q^{2n+2}}{cdef}} \\ = 
  \frac{(aq, \frac qa, \frac{aq}{cd}, \frac{aq}{ef};q)_n}
  { (\frac qc, \frac qd, \frac{aq}{e}, \frac{aq}{f};q)_n} 
  \bbhs{4}{4} {e,f, \frac{aq^{n+1}}{cd}, q^{-n}}{ \frac{aq}{c}, \frac{aq}{d}, q^{n+1}, \frac{ef}{aq^n}  }
  {q}{q}. \label{8psi8trans}
\end{multline}

Note further that~\eqref{8psi8trans} is a bilateral analog of Watson's $q$-analog of
Whipple's theorem~\cite{W29}(cf.~\cite[p. 360, Eq. (III.17)]{GR04})

Notice that for level $\ell = 1$ and $2$, $q^{-n}(q)_n(q^2;q)_n \beta^{(\ell,1)}_n$ is a
limiting case of a ${}_{10}\psi_{10}$, while for $\ell>6$, it is a limiting case of a ${}_{t}\psi_{t}$
with $t = \ell - 1 + (1+(-1)^\ell)/2.$

More precisely, if $\ell$ is even and $\ell\geqq 6$,   
\begin{multline*} q^{-n } (q)_n (q^2;q)_{n+1} \beta^{(\ell,1)}_n(1,q) \\ =
  \displaystyle{ \underset{e\to \infty}{\lim} \bbhs{\ell}{\ell}
  {q^{\frac 72},-q^{\frac 72},q^{-n},q^{1-n},q^{2-n},\overbrace{e,e,\dots,e}^{\ell-6},-q^2}
  {q^{\frac 12},-q^{\frac 12},q^{n+4},q^{n+3},q^{n+2},\underbrace{\frac{q^4}{e},\frac{q^4}{e},\dots,\frac{q^4}{e}}_{\ell-6},-q^2}{q^3}
  { \frac{-q^{3n+2\ell-12}}{e^{\ell-6}} } },
\end{multline*}
while if $\ell$ is odd and $\ell>6$,
\begin{multline*}
q^{-n } (q)_n (q^2;q)_{n+1} \beta^{(\ell,1)}_n(1,q) \\ =
  \displaystyle{ \underset{e\to \infty}{\lim} \bbhs{\ell-1}{\ell-1}
  {q^{\frac 72},-q^{\frac 72},q^{-n},q^{1-n},q^{2-n},\overbrace{e,e,\dots,e}^{\ell-6}}
  {q^{\frac 12},-q^{\frac 12},q^{n+4},q^{n+3},q^{n+2},\underbrace{\frac{q^4}{e},\frac{q^4}{e},\dots,\frac{q^4}{e}}_{\ell-6}}{q^3}
  { \frac{q^{3n+2\ell-12}}{e^{\ell-6}} } }.
\end{multline*}

Then, to obtain the series and product expressions for $\chi( \Omega(\ell\Lambda_0))$, one inserts
the Bailey pair $\Big(\alpha^{(\ell,1)}_n(1,q), \beta^{(\ell,1)}_n(1,q)\Big)$ into~\eqref{WBL} with $a=1$,
 and upon applying~\eqref{JTP} and~\eqref{QPI}, we find that
 \begin{multline} \label{GenSerProdL0}
   \sum_{m=0}^\infty q^{9m^2} \Big(q^{-6m+1}  \beta^{(\ell,1)}_{3m-1}(1,q)+ \beta^{(\ell,1)}_{3m}(1,q) + q^{6m+1}  \beta^{(\ell,1)}_{3m+1}(1,q)
   \Big) \\
   = \frac{(q,q^{\ell+2}, q^{\ell+3}; q^{\ell+3})_\infty ( q^{\ell+1}, q^{ \ell+ 5} ; q^{2\ell+6}) }{(q)_\infty}.
   \end{multline}
   
   And thus in~\eqref{GenSerProdL0}, 
 we have a uniform series-product identity for the principally specialized character
 of the $(\ell, 0)$ standard module of $A_2^{(2)}$ for any $\ell$.
 
  To express the $\beta_n^{(\ell,1)}$ as a multisum for arbitrary $\ell$, 
one may employ the Andrews--Baxter--Forrester
bilaterial very-well poised $q$-hypergeometric summation formula~\cite[p. 83, Eq. (8.56)]{A86};
cf.~\cite[Appendix B, pp. 261--265]{ABF84}. 

  The level 3 case will be considered in detail in the next section.

\subsection{Bailey pairs for $\chi (\Omega( (\ell-2)\Lambda_0 + \Lambda_1 )) $}
 \[ \alpha_n = \alpha_n^{(\ell,2)}(1,q) = \left\{ 
\begin{array}{ll} 1 & \mbox{if $n=0$}\\
q^{\frac 32(\ell-3)r^2 - \frac12 (9-\ell) r} +q^{\frac 32(\ell-3)r^2 + \frac 12(9-\ell)r } & \mbox{if $n=3r>0$ } \\
-q^{\frac 32(\ell-3)r^2 + \frac12 (\ell+3) r + 1} &\mbox{if $n=3r+1$} \\
-q^{\frac 32(\ell-3)r^2 - \frac12 (\ell+3) r + 1} &\mbox{if $n=3r-1$}
\end{array}
 \right. .
\]

The calculation of $\beta_n^{(\ell,2)}$ parallels that of $\beta_n^{(\ell,1)}$.  The details of
the $\ell=7$ case are given by Slater~\cite[p. 464]{S51}.
\begin{multline}
\beta_{n}^{(\ell,2)} (1,q)  = \sum_{s=0}^n \frac{\alpha_s^{(\ell,2)} (1,q) }{(q)_{n-s} (q)_{n+s}}     \\
    =  \frac{1}{ (q)_n (q^2;q)_{n} } 
                  \sum_{r\in\Z}\frac{  
                  ( 1- q^{6r+1} ) (q^{-n};q)_{3r}   }{ (1-q) (q^{n+2};q)_{3r} } (-1)^r q^{\frac 32(\ell-6)r^2 + \frac 12(\ell-6)r } .
\end{multline}
Notice that 
\begin{equation} \label{beta12}
q^{n} \beta_n^{(\ell,2)}(1,q) = \beta_n^{(\ell,1)}(1,q). 
\end{equation}

And so it follows that the series and product expressions for $\chi(\Omega((\ell-2)\Lambda_0 + \Lambda_1))$ are
 \begin{multline}  \label{GenSerProdL1}
   \sum_{m=0}^\infty q^{9m^2} \Big(q^{-6m+1}  \beta^{(\ell,2)}_{3m-1}(1,q)+ \beta^{(\ell,2)}_{3m}(1,q) + q^{6m+1}  \beta^{(\ell,2)}_{3m+1}(1,q)
   \Big) \\
   = \frac{(q^2,q^{\ell+1}, q^{\ell+3}; q^{\ell+3})_\infty ( q^{\ell-1}, q^{ \ell+ 7} ; q^{2\ell+6}) }{(q)_\infty},
   \end{multline}
 a general identity corresponding to the $(\ell-2, 1)$ module.

\section{Level $3$}
Let us consider the case $\ell=3$ in detail. 
This level is of particular interest as it was the study of the the level $3$ standard modules
of $A_2^{(2)}$ that led S. Capparelli to discover
two new 
Rogers--Ramanujan type partition identities~\cite{AAG95, A94,C88,C96}.  See~\cite[\S 3]{S10} for some
historical notes.

From the $(3,0)$-module, Capparelli conjectured (and later proved~\cite{C96}, although the first proof was due to
Andrews~\cite{A94}) the following partition identity.   A \emph{partition} $\lambda$ of an integer $n$ is a finite
weakly decreasing sequence $(\lambda_1, \lambda_2, \dots, \lambda_l)$ of positive integers that sum to $n$; 
each $\lambda_i$ is called a \emph{part} of the partition $\lambda$.
\begin{theorem}[Capparelli's first partition identity]

Let $c_1(n)$ denote the number of partitions $\lambda = (\lambda_1, \dots, \lambda_l)$ of $n$ wherein 
   \begin{itemize}
     \item  $\lambda_i \neq 1$ for $i=1,2,\dots,l$,
     \item $\lambda_i - \lambda_{i+1}\geq 2$, for $i=1,2,\dots, l-1$,
     \item $\lambda_i - \lambda_{i+1} =2$ only if $\lambda_i \equiv 1\pmod{3}$,
     \item $\lambda_i - \lambda_{i+1} =3 $ only if $\lambda_i \equiv 0 \pmod{3}$.
   \end{itemize}
   
Let $c_2(n)$ denote the number of partitions of $n$ into distinct parts $\not\equiv\pm 1\pmod{6}$.  
 Let $c_3(n)$ denote the number of partitions of $n$ into parts congruent to $\pm 2,\pm 3 \pmod{12}$.
 Then $c_1(n) = c_2(n) = c_3(n)$ for all $n$.
\end{theorem}

From the $(1,1)$-module, Capparelli obtained the companion identity:
\begin{theorem}[Capparelli's second partition identity]

Let $d_1(n)$ denote the number of partitions $\lambda = (\lambda_1, \dots, \lambda_l)$ of $n$ wherein 
   \begin{itemize}
     \item  $\lambda_i \neq 2$ for $i=1,2,\dots,l$,
     \item $\lambda_i - \lambda_{i+1}\geq 2$, for $i=1,2,\dots, l-1$,
     \item $\lambda_i - \lambda_{i+1} =2$ only if $\lambda_i \equiv 1\pmod{3}$,
     \item $\lambda_i - \lambda_{i+1} =3 $ only if $\lambda_i \equiv 0 \pmod{3}$.
   \end{itemize}
Let $d_2(n)$ denote the number of partitions of $n$ into distinct parts $\not\equiv\pm 2\pmod{6}$.  
   
 Then $d_1(n) = d_2(n)$ for all $n$.
\end{theorem}

\subsection{ $\chi (\Omega( 3\Lambda_0 ))$ }
 In order to use~\eqref{8psi8trans}, we need to 
consider three cases, $n= 3m$, $3m+1$, and $3m-1$.

In~\eqref{8psi8trans}, replace $q$ by $q^3$; then set $a=q$, $n=m$, $f=q^{2-3m}$, and $c=d=e$,
to obtain

\begin{align}
\beta^{(3,1)}_{3m}(1,q) &= \frac{q^n(1-q)}{(q)_n (q)_{n+1}} \lim_{e\to 0} 
\bbhs{8}{8}{ q^{-7/2}, -q^{7/2}, e,e,e,q^{2-3m}, q^{1-3m},q^{-3m}}{q^{1/2},-q^{1/2}, 
\frac{q^4}{e}, \frac{q^4}{e}, \frac{q^4}{e}, q^{3m+2}, q^{3m+3}, q^{3m+4}  }{q^3}{q^3} \notag \\
&= \frac{q^{3m}(1-q)}{(q)_{3m} (q)_{3m+1}}\lim_{e\to 0}\frac{(q^4,q^2, \frac{q^4}{e^2}, \frac{q^{3m+2}}{e};q^3)_m}
{(  \frac{q^3}{e}, \frac{q^3}{e}, \frac{q^4}{e}, q^{3m+2}   ;q^3)_m}
\bbhs{4}{4}{e, q^{2-3m}, \frac{q^{3m+4}}{e^2}, q^{-3m}}{ \frac{q^4}{e}, \frac{q^4}{e}, q^{3m+3}, eq^{1-6m}}
{q^3}{q^3} \notag\\
&= \sum_{r=-m}^m \frac{ (-1)^{m+r} q^{\frac 32 m^2 - \frac 12 m - \frac 12 r - 3mr + \frac 32 r^2} }
{ (q^2;q^3)_{2m} (q^3;q^3)_{m+r} (q;q^3)_{m-r} (q^3;q^3)_{m-r}  } \notag\\
& = \sum_{r=0}^{2m} \frac{(-1)^r q^{\frac 32 r^2 + \frac 12 r} (q^2;q^3)_r}{(q^2;q^3)_{2m} (q^3;q^3)_{2m-r} (q)_{3r}}. \label{CapBeta3m}
\end{align}

In~\eqref{8psi8trans}, replace $q$ by $q^3$; then set $a=q$, $n=m$, $f=q^{-1-3m}$, and $c=d=e$,
to obtain

\begin{align}
\beta^{(3,1)}_{3m+1}(1,q) &= \frac{q^n(1-q)}{(q)_n (q)_{n+1}} \lim_{e\to 0} 
\bbhs{8}{8}{ q^{-7/2}, -q^{7/2}, e,e,e,q^{-1-3m}, q^{1-3m},q^{-3m}}{q^{1/2},-q^{1/2}, 
\frac{q^4}{e}, \frac{q^4}{e}, \frac{q^4}{e}, q^{3m+2}, q^{3m+5}, q^{3m+4}  }{q^3}{q^3}  \notag\\
&= \frac{q^{3m+1}(1-q)}{(q)_{3m+1} (q)_{3m+2}}\lim_{e\to 0}\frac{(q^4,q^2, \frac{q^4}{e^2}, \frac{q^{3m+5}}{e};q^3)_m}
{(  \frac{q^3}{e}, \frac{q^3}{e}, \frac{q^4}{e}, q^{3m+5}   ;q^3)_m}
\bbhs{4}{4}{e, q^{-1-3m}, \frac{q^{3m+4}}{e^2}, q^{-3m}}
{ \frac{q^4}{e}, \frac{q^4}{e}, q^{3m+3}, eq^{-2-6m}}{q^3}{q^3} \notag\\
&= \sum_{r=-m}^m \frac{ (-1)^{m+r} q^{\frac 32 m^2 + \frac 72 m - \frac 72 r - 3mr + \frac 32 r^2 + 1} }
{ (q^2;q^3)_{2m+1} (q^3;q^3)_{m+r} (q;q^3)_{m-r+1} (q^3;q^3)_{m-r}  }\notag\\
& = \sum_{r=0}^{2m} \frac{(-1)^r q^{\frac 32 r^2 + \frac 72 r +1}(q^2;q^3)_r}{(q^2;q^3)_{2m+1} (q^3;q^3)_{2m-r}(q)_{3r+1}}. \label{CapBeta3mPlus1}
\end{align}

For convenience, let us define the abbreviation
 \begin{equation}
   \sigma(m,r):= \frac{(-1)^r q^{\frac 32 r^2 + \frac 12 r} (q^2;q^3)_r }{ (q^2;q^3)_{2m} (q^3;q^3)_{2m-r}
    (q)_{3r} },
 \end{equation}
so that we have immediately
   \[ \beta^{(3,1)}_{3m}(1,q) = \sum_{r=0}^{2m} \sigma(m,r), \]  and with a bit of elementary algebra,
   \[ \beta^{(3,1)}_{3m+1}(1,q) = \sum_{r=0}^{2m} \frac{\sigma(m,r)}{1-q^{6m+2}} \left( \frac{1}
   {1-q^{3r+1}} - 1 \right), \]
   for $m\geqq 0$.

The author could not find a direct substitution into ~\eqref{8psi8trans}, analogous
to the $n=3m$ and $n=3m+1$ cases, which yields
the $n=3m-1$ case.  So we resort to an alternate method to obtain the $n=3m-1$ case.

From the Paule--Riese \texttt{qZeil.m} \emph{Mathematica} package available for download at \texttt{http://www.risc.jku.at/research/combinat/software/qZeil/index.php} and documented
in~\cite{PR97}, one can find that $\beta^{(3,1)}_n(1,q)$ satisfies the recurrence
\begin{equation} \label{CapBetaRecN}
\beta_n = \frac{-q^2 + q^{2n} + q^{2n+1}}{q^2(1-q^{2n})(1-q^{2n-1}) }\beta_{n-1} - \frac{1}{(1-q^{2n})(1-q^{2n-1})} \beta_{n-2}
\end{equation}
as certified by the rational function
\[ \frac{q^{-n-6 r-2} \left(q^{3 r}-q^n\right) \left(q^{3 r+1}-q^n\right) \left(q^{3
   r+2}-q^n\right)}{\left(q^n-1\right) \left(q^n+1\right) \left(q^{2 n}-q\right)
   \left(q^{6 r+1}-1\right)} . \]



Setting $n=3m+1$ in~\eqref{CapBetaRecN} and rearranging, we see how to express 
$\beta_{3m-1}$ in terms of the two known expressions $\beta_{3m}$ and $\beta_{3m+1}$:
\begin{align} \label{CapBetaRecM}
\beta_{3m-1} &= -(1-q^{6m+2})(1-q^{6m+1})\beta_{3m+1} - (1 - q^{6m} - q^{6m+1})\beta_{3m}\\
 &=  -(1-q^{6m+2})(1-q^{6m+1}) \sum_{r=0}^{2m} \frac{\sigma(m,r)}{1-q^{6m+2}} \left( \frac{1}
   {1-q^{3r+1}} - 1 \right) \notag\\ & \qquad\qquad+ (1-q^{6m+1})\beta_{3m} - (1-q^{6m}-q^{6m+1}) \beta_{3m} \notag\\
 & = -(1-q^{6m+1}) \sum_{r=0}^{2m} \frac{\sigma(m,r)}{1-q^{3r+1}} + q^{6m}\beta_{3m} \notag\\
 & = \sum_{r=0}^{2m} \sigma(m,r) \left( q^{6m} - \frac{1-q^{6m+1}}{1-q^{3r+1}} \right) ,
 \label{CapBeta3mMinus1}
 \end{align} for $m\geqq 1$.
 
 Inserting $\Big(\alpha^{(3,1)}_n(1,q), \beta^{(3,1)}_n(1,q)\Big)$ into~\eqref{WBL} with $a=1$,
 and applying~\eqref{JTP} and~\eqref{QPI}, we find that
 \begin{multline} \label{Cap1}
   \sum_{m=0}^\infty q^{9m^2} \Big(q^{-6m+1}  \beta^{(3,1)}_{3m-1}(1,q)+ \beta^{(3,1)}_{3m}(1,q) + q^{6m+1}  \beta^{(3,1)}_{3m+1}(1,q)
   \Big) \\
 = \sum_{m=0}^\infty \sum_{r=0}^{2m} q^{9m^2} \sigma(m,r) 
  \left( 
  q^{1-6m}\left( q^{6m} - \frac{1-q^{6m+1}}{1-q^{3r+1}} \right)
  + 1 \right. \\
 \left.  + \frac{q^{6m+1}}{1-q^{6m+2}}\left( 1 - \frac{1}{1-q^{3r+1}} \right)
  \right) \\
   =\frac{(q,q^5,q^6;q^6)_\infty(q^4,q^8;q^{12})_\infty}{(q)_\infty}= \frac{1}{(q^2,q^3,q^9,q^{10};q^{12})}_\infty \\
   = (-q^2;q^2)_\infty (-q^3;q^6)_\infty.
 \end{multline}
 The series expansion of $(q^2,q^3,q^9,q^{10};q^{12})_\infty^{-1}$ in~\eqref{Cap1} is quite
 different than others that have appeared in the 
 literature, due to Alladi, Andrews, and Gordon~\cite[pp. 648--649, Lemma 2(b)]{AAG95}
 (cf. \cite[p. 399, Eq. (1.3)]{S04}), the author~\cite[p. 399, Eq. (1.4) and Eq. (1.5)]{S04},
 and Bringmann and Mahlburg~\cite{BM15}.
 
 \subsection{$ \chi (\Omega( \Lambda_0 +\Lambda_1 ))$  }

  In light of~\eqref{beta12}, it is trivial to obtain $\beta^{(3,2)}_n(1,q)$ from $\beta^{(3,1)}_n(1,q)$,
and upon inserting $\Big(\alpha^{(3,2)}_n(1,q), \beta^{(3,2)}_n(1,q)\Big)$ into~\eqref{WBL} with $a=1$,
 and applying~\eqref{JTP} and~\eqref{QPI}, we find that
 \begin{multline} \label{Cap2}
  \sum_{m=0}^\infty \sum_{r=0}^{2m} q^{9m^2-3m} \sigma(m,r) 
  \left( 
  q^{2-6m}\left( q^{6m} - \frac{1-q^{6m+1}}{1-q^{3r+1}} \right)
  + 1 \right. \\
 \left. + \frac{q^{6m}}{1-q^{6m+2}}\left( 1 - \frac{1}{1-q^{3r+1}} \right)
  \right) \\
   =\frac{(q^2,q^4,q^6;q^6)_\infty(q^2,q^{10};q^{12})_\infty}{(q)_\infty}
   =(-q;q^2)_\infty (-q^6;q^6)_\infty.
 \end{multline}
 
 \subsection{Nandi's recent work on level $4$}
 It should be noted that recently D. Nandi, in his Ph.D. thesis~\cite{N14} conjectured 
the partition identities corresponding to the three inequivalent level $4$ standard modules
$(4,0)$, $(2,1)$ and $(0,2)$.  These identities, while still in the spirit of the Rogers--Ramanujan
and Capparelli identities, involve difference conditions that are \emph{much} more 
complicated than anything that has been considered previously in the theory of
partitions.  It is no wonder that after Capperelli's discoveries for level $3$, it took a 
quarter century to successfully perform the analogous feat for level $4$.
    
\section{Bailey pairs for levels $3$ through $9$ summarized }
\subsection{Level 3}
\begin{align*}
\beta^{(3,1)}_{3m}(1,q) &= \sum_{r=-m}^m \frac{(-1)^{m+r} q^{\frac 32 m^2 - \frac 12 m +3mr + \frac 32 r^2 -
\frac 12 r} (q^2;q^3)_{2m} }{ (q)_{6m} (q^2;q^3)_{m+r}   } \qbin{2m}{m+r}{q^3}\\
&=\sum_{r=0}^{2m} \frac{(-1)^r q^{\frac 32r^2 - \frac 12 r}}{(q;q^3)_{2m} (q^2;q^3)_r (q^3;q^3)_{2m-r}
(q^3;q^3)_r }
\end{align*}
\begin{align*}
\beta^{(3,1)}_{3m+1}(1,q) &= \sum_{r=-m}^m \frac{(-1)^{m+r} q^{\frac 32 m^2 + \frac 52 m +3mr + \frac 32 r^2 +
\frac 52 r + 1} (q^2;q^3)_{2m} }{ (q)_{6m+1} (q^2;q^3)_{m+r}   } \qbin{2m}{m+r}{q^3}\\
&=\sum_{r=0}^{2m} \frac{(-1)^r q^{\frac 32r^2 + \frac 52 r}}{(q;q^3)_{2m+1} (q^2;q^3)_{r+1} (q^3;q^3)_{2m-r}
(q^3;q^3)_r }
\end{align*}
\begin{equation*}
\beta_{3m-1} = -(1-q^{6m+1})(1-q^{6m+2}) \beta_{3m+1} - (1-q^{6m}-q^{6m+1})\beta_{3m}
\end{equation*}

\subsection{Level 4}
\begin{align*}
\beta^{(4,1)}_{3m}(1,q) &= \sum_{r=-m}^m \frac{(-1)^{m+r} q^{3 m^2 +3mr + \frac 32 r^2 -
\frac 12 r} (-q^2;q^3)_r (q^2;q^3)_{2m} }{ (q)_{6m} (-q^2;q^3)_m (q^2;q^3)_{m+r}   } \qbin{2m}{m+r}{q^3}\\
&=\sum_{r=0}^{2m} \frac{(-1)^r q^{3m^2 -3mr + 3r^2} }
{(-q;q^3)_{m-r} (q;q^3)_{2m} (-q^2;q^3)_m (q^2;q^3)_r (q^3;q^3)_{2m-r}
(q^3;q^3)_r }
\end{align*}

\begin{align*}
\beta^{(4,1)}_{3m+1}(1,q) &= \sum_{r=-m}^m \frac{(-1)^{m+r} q^{3 m^2 +3m  +3mr + \frac 32 r^2 +
\frac 52 r + 1} (-q^2;q^3)_r (q^2;q^3)_{2m} }{ (q)_{6m+1} (-q^2;q^3)_m (q^2;q^3)_{m+r+1}   } \qbin{2m}{m+r}{q^3}\\
&=\sum_{r=0}^{2m} \frac{(-1)^r q^{3m^2 - 3mr + 3r^2 + 3r + 1}  }
{(-q;q^3)_{m-r} (q;q^3)_{2m+1}(-q^2;q^3)_m  (q^2;q^3)_{r+1} (q^3;q^3)_{2m-r}
(q^3;q^3)_r }
\end{align*}

\begin{equation*}
q^{3m} \beta_{3m-1} = -(1-q^{6m+1})(1-q^{6m+2})(1+q^{3m}) \beta_{3m+1} -q^{3m} (1+q^{3m})
(1-q^{3m}-q^{3m+1})\beta_{3m}
\end{equation*}

\subsection{Level 5}
\begin{equation*}
\beta^{(5,1)}_{3m}(1,q)  = \frac{ q^{9m^2} } {(q)_{6m} },
\end{equation*}

\begin{equation*}
\beta^{(5,1)}_{3m+1}(1,q) = \frac{ q^{9m^2 + 6m + 1}}{ (q)_{6m+2} },
\end{equation*} and

\begin{equation*}
q^{6m-1} \beta_{3m-1} = (1-q^{6m})(1-q^{6m-1}) \beta_{3m} ,
\end{equation*}
which, together, simplifies to Eq.~\eqref{beta51}.

\subsection{Level 6}
\begin{equation*}
\beta^{(6,1)}_{3m}(1,q)  = \frac{ q^{3m}  (-1;q)_{3m} } {(-1;q^3)_{3m} (q)_{6m} },
\end{equation*}

\begin{equation*}
\beta^{(6,1)}_{3m+1}(1,q) = \frac{ q^{3m + 1} (-1;q^3)_{3m+1} }{ (q)_{6m+2} (-1;q)_{3m+1}},
\end{equation*}
and
\begin{equation*}
(q+q^{6m-1}-q^{3m})\beta_{3m-1} = (1-q^{6m})(1-q^{6m-1}) \beta_{3m} , 
\end{equation*}
and thus~\eqref{beta61} holds.

\subsection{Level 7}
\begin{equation*}
\beta^{(7,1)}_{3m}(1,q)  = \frac{ q^{3m} } {(q)_{6m} },
\end{equation*}

\begin{equation*}
\beta^{(7,1)}_{3m+1}(1,q) = \frac{ q^{3m+ 1}}{ (q)_{6m+2} },
\end{equation*}
and
\begin{equation*}
q\beta_{3m-1} = (1-q^{6m})(1-q^{6m-1})\beta_{3m},
\end{equation*}
and thus~\eqref{beta71} holds.

\subsection{Level 8}

\begin{align*}
\beta^{(8,1)}_{3m}(1,q) &= \sum_{r=-m}^m \frac{ q^{ \frac 32 r^2 + \frac 12 r+3m} (-q^2;q^3)_r (q^2;q^3)_{2m} }{ (q)_{6m} (-q^2;q^3)_m (q^2;q^3)_{m+r}   } 
\qbin{2m}{m+r}{q^3}\\
&=\sum_{r=0}^{2m} \frac{q^{ 3m^2 + 2m + 3r^2 + r - 6mr}  }
{(-q;q^3)_{m-r} (q;q^3)_{2m} (-q^2;q^3)_m (q^2;q^3)_r (q^3;q^3)_{2m-r} (q^3;q^3)_r },
\end{align*}

\begin{align*}
\beta^{(8,1)}_{3m+1}(1,q) &= \sum_{r=-m}^m \frac{q^{ \frac 32 r^2 +
\frac 12 r +3m+ 1} (-q^2;q^3)_r (q^2;q^3)_{2m} }{ (q)_{6m+1} (-q^2;q^3)_m (q^2;q^3)_{m+r+1}   } \qbin{2m}{m+r}{q^3}\\
&=\sum_{r=0}^{2m} \frac{q^{ 3m^2+2m+3r^2+r-6mr+1 }  }
{(-q;q^3)_{m-r} (q;q^3)_{2m+1}(-q^2;q^3)_m  (q^2;q^3)_{r+1} (q^3;q^3)_{2m-r} (q^3;q^3)_r },
\end{align*}

\begin{equation*}
q^{2} \beta_{3m-1} = -(1-q^{6m+1})(1-q^{6m+2})(1+q^{3m}) \beta_{3m+1} + q
(1+q+q^{3m}-q^{6m+1})\beta_{3m}.
\end{equation*}

\subsection{Level 9}
\begin{align*}
\beta^{(9,1)}_{3m}(1,q) &= \sum_{r=-m}^m \frac{ q^{ 3 r^2 + r+3m}  (q^2;q^3)_{2m} }
{ (q)_{6m} (q^2;q^3)_{m+r} }
\qbin{2m}{m+r}{q^3}\\
&=\sum_{r=0}^{2m} \frac{q^{ 3m^2+2m+3r^2+r-6mr} }
{(q;q^3)_{2m} (q^2;q^3)_r (q^3;q^3)_{2m-r}
(q^3;q^3)_r },
\end{align*}
\begin{align*}
\beta^{(9,1)}_{3m+1}(1,q) &= \sum_{r=-m}^m \frac{ q^{ 3 r^2 + r+3m +1}  (q^2;q^3)_{2m} }
{ (1-q^{3m+1})(1-q^{3m+2})(q)_{6m} (q^2;q^3)_{m+r} }
\qbin{2m}{m+r}{q^3}\\
&=\sum_{r=0}^{2m} \frac{q^{ 3m^2+2m+3r^2+r-6mr+1} }
{(1-q^{3m+1})(1-q^{3m+2})(q;q^3)_{2m} (q^2;q^3)_r (q^3;q^3)_{2m-r}
(q^3;q^3)_r },
\end{align*}\begin{equation*}
q^{3} \beta_{3m-1} = -(1-q^{6m+1})(1-q^{6m+2}) \beta_{3m+1} + q
(1+q-q^{6m+1})\beta_{3m}.
\end{equation*}

\section{Conclusion and Open Questions}
It is the hope of the author that the results presented here will help to provide some 
insight into the structure of $A_2^{(2)}$ that can be exploited by vertex operator
algebraists.  Questions of course remain.  For instance, 
as pointed out by Ole Warnaar during the question-and-answer period following my
talk at the Alladi conference, 
it is not at all clear 
how the series expressions in~\eqref{Cap1} and~\eqref{Cap2} enumerate the 
partition functions $c_1(n)$ and $c_2(n)$ respectively. 
It would be very nice indeed if this connection could be established.  
Christian Krattenthaler pointed out that it was conceivable that there are 
\emph{other} families of Bailey pairs that could give rise to identities with
the same product sides.
While this is true, the choice of the $\alpha_n$ employed here is motivated by
classical work; in particular the level 5 and 7 identities and the corresponding Bailey pairs
coincide with the work of Slater~\cite{S51,S52}.  Further the availability of the 
Andrews--Baxter--Forrester transformation to express the $\beta_n$ as a
multisum for \emph{any} level $\ell$ is an encouraging sign that this may be 
a fruitful direction to pursue in the effort to better understand $A_2^{(2)}$ as a whole.
 
\section*{Acknowledgments}
The author thanks Jim Lepowsky and Robert Wilson for assistance with the
exposition in Section 1.2.


\begin{thebibliography}{29}

\bibitem{AAG95} K. Alladi, G. E. Andrews, and B. Gordon, Refinements and generalizations
of Capparelli's conjecture on partitions, J. Algebra 174 (1995) 636--658.

\bibitem{AAB87} A. K. Agarwal, G. E. Andrews, and D. M. Bressoud,
The Bailey lattice, J. Indian Math. Soc. (N. S.) 51 (1987) 57--73.

\bibitem{A74} G. E. Andrews,
An analytic
generalization of the Rogers--Ramanujan
identities for odd moduli,
Proc. Nat. Acad. Sci. USA,
 {71} (1974) 4082--4085.


\bibitem{A84} G. E. Andrews,
{Multiple series Rogers--Ramanujan
type identities},
Pacific J. Math.,  {114} (1984)
267--283.

\bibitem{A86} G. E. Andrews,
\emph{q-series: their development and
application in analysis, number theory, combinatorics, physics, and
computer algebra}, CBMS Regional Conferences Series in Mathematics, no. 66,
American Mathematical Society, Providence, RI, 1986.

\bibitem{A94} G. E. Andrews, Schur's theorem,
Capparelli's conjecture and $q$-trinomial coefficients.
The Rademacher legacy to mathematics (University Park, PA, 1992) 141--154,
Contemp. Math., 166, Amer. Math. Soc., Providence, RI, 1994.

\bibitem{ABF84} G. E. Andrews, R. J. Baxter, P. J. Forrester,
Eight-vertex SOS model and generalized Rogers--Ramanujan type identities,
J. Stat. Phys. 35 (1984) 193--266.

\bibitem{B36} W. N. Bailey, Series of hypergeometric type which are infinite in both 
directions, Quart. J. Math. (Oxford) 7 (1936) 105--115.

\bibitem{B49} W.~N.~Bailey, {Identities of the Rogers-Ramanujan type},
Proc. London Math. Soc. (2) {50} (1949) 1--10.

\bibitem{B79} D. M. Bressoud,
A generalization of the Rogers-Ramanujan identities for all moduli,
J. Combin. Theory Ser. A 27 (1979) 64--68.

\bibitem{B80} D. M. Bressoud,
Analytic and combinatorial generalizations of the Rogers-Ramanujan
identities, Mem. Amer. Math. Soc. 24 (1980) no. 227, 1--54.

\bibitem{BM15} K. Bringmann and K. Mahlburg,
False theta functions and companions to Capparelli's identities,
Adv. Math. 278 (2015) 121--136.

\bibitem{C88} S. Capparelli, Vertex operator relations for affine algebras and combinatorial
identities, Ph.D thesis, Rutgers University, 1988.

\bibitem{C96} S. Capparelli, A construction of the level $3$ modules for the
affine Lie algebra $A^{(2)}_2$ and a new combinatorial identity of the
Rogers-Ramanujan type, Trans. Amer. Math. Soc. 348 (1996) 481--501.

\bibitem{F16} R. Fricke, Die Elliptischen Funktionen und ihre Anwendungen,
Ers Teil, Teubner, Leipzig, 1916.

\bibitem{GR04} G. Gasper and M. Rahman, \emph{Basic Hypergeometric
Series}, 2nd ed., Cambridge University Press, 2004.

\bibitem{G61} B. Gordon,  A combinatorial generalization of the
Rogers-Ramanujan identities, Amer. J. Math. 83 (1961) 393--399.

\bibitem{K90} V. G. Kac, \emph{Infinite Dimensional Lie Algebras}, 3rd ed.,
Cambridge Univ. Press, 1990.

\bibitem{KKLW81} V. G. Kac, D. A. Kazhdan, J. Lepowsky, and R. Wilson,
Realization of the basic representations of the Euclidean Lie algebras,
Adv. Math. 42 (1981) 83--112.


\bibitem{LM78} J. Lepowsky and S. Milne,
Lie algebraic approaches to classical partition identities,
Adv. Math. 29 (1978) 15--59.

\bibitem{LW78} J. Lepowsky and R. L. Wilson,
Construction of the affine Lie algebra $A_1^{(1)}$,
Comm. Math. Phys. 62 (1978) 43--53.

\bibitem{LW82} J. Lepowsky and R. L. Wilson,
A Lie theoretic interpretation and proof of the Rogers-Ramanujan
identities, Adv. Math. 45 (1982) 21--72.

\bibitem{LW84} J. Lepowsky and R. L. Wilson,
The structure of standard modules I: Universal algebras and the
Rogers-Ramanujan identities, Invent. Math. 77 (1984) 199--290.

\bibitem{LW85} J. Lepowsky and R. L. Wilson,
The structure of standard modules II: the case $A_1^{(1)}$, principal
gradation, Invent. Math. 79 (1985) 417--442.

\bibitem{MS08} J. McLaughlin and A. V. Sills,
Ramanujan--Slater type identities related to the moduli 18 and 24,
J. Math. Anal. Appl. (2008) 765--777.

\bibitem{N14} D. Nandi, Partition identities arising from the standard $A_2^{(2)}$
modules of level 4, Ph.D. thesis, Rutgers University, 2014.

\bibitem{PR97} P. Paule and A. Riese, A Mathematica $q$-analogue of
Zeilberger's algorithm based on an algebraically motivated approach to
$q$-hypergeometric telescoping, in \emph{Special Functions, $q$-Series
and Related Topics}, Fields Inst. Commun., vol. 14 (1997) 179--210.

\bibitem{S04} A. V. Sills, On series expansions of Capparelli's infinite product,
Adv. Appl. Math. 33 (2004) 397--408.

\bibitem{S10} A. V. Sills, Rademacher-type formulas for restricted partition
and overpartition functions, Ramanujan J. 23 (2010) 253--264.

\bibitem{S51} L. J. Slater, {A new proof of Rogers transformation of
infinite series}, Proc. London Math Soc. (2) {53} (1951)
460--475

\bibitem{S52} L. J. Slater, {Further identities of the Rogers-Ramanujan
type}, Proc. London Math Soc. (2) {54} (1952) 147--167.

\bibitem{W29}  G. N. Watson, {A new proof of the Rogers--Ramanujan identities},
J. London Math. Soc. 4 (1929) 4--9.


\end{thebibliography}
\end{document}